\documentclass[12pt]{article}
\usepackage{mathrsfs}
\usepackage{amsmath, amsfonts, amscd, amssymb}
\usepackage{epsfig, graphicx, texdraw, color, latexsym, colordvi, cite}

\newtheorem{theo}{Theorem}[section]

\newtheorem{exam}[theo]{ Example}

\def\qed{\hfill \rule{4pt}{7pt}}
\def\pf{\noindent {\it Proof.} }

\numberwithin{equation}{section}

\def\x{{\mathbf x}}
\def\y{{\mathbf y}}

\def\pa{{\rm Par}}

\begin{document}
\begin{center}
{\bf \Large Partially Ordinal Sums and $P$-partitions}

\vskip 20pt
Daniel K. Du\raisebox{5pt}{\scriptsize 1} and
Qing-Hu Hou\raisebox{5pt}{\scriptsize 2}

\vskip 20pt
Center for Combinatorics, LPMC-TJKLC \\
Nankai University, Tianjin 300071, P. R. China

\vskip 3mm

\raisebox{5pt}{\scriptsize 1}dukang@mail.nankai.edu.cn,
\raisebox{5pt}{\scriptsize 2}hou@nankai.edu.cn,

\end{center}

\begin{abstract}
We present a method of computing the generating function $f_P(\x)$ of $P$-partitions of a poset $P$.
The idea is to introduce two kinds of transformations on posets and compute $f_P(\x)$ by recursively applying these transformations. As an application, we consider the partially ordinal sum $P_n$ of $n$ copies of a given poset, which generalizes both the direct sum and the ordinal sum. We show that the sequence $\{f_{P_n}(\x)\}_{n\ge 1}$ satisfies a finite system of recurrence relations with respect to $n$. We illustrate the method by several examples, including a kind of $3$-rowed posets and the multi-cube posets.

\end{abstract}



\section{Introduction}
A $P$-partition is an order-reversing map from a partially ordered set (poset) to
non-negative integers (see \cite[Ch. IV]{Stanley-97}). Denote the set of $P$-partitions of a poset $P$ by $\pa(P)$. The (multivariate) generating function of $P$-partitions of $P$ is
given by
\[
f_P(\x) = \sum_{\sigma \in \pa(P)} \prod_{a \in P} x_a^{\sigma(a)},
\]
where $\x$ is the variable vector $(x_a)_{a\in P}$.

Stanley showed that $f_P(\x)$ can be expressed as a sum over linear extensions of $P$ \cite[Theorem 4.5.4]{Stanley-97}. Andrews, Paule and Riese \cite{Andre01-a, Andre01-c,
Andre01-e, Andre04, Andre04-a,Andre07b}
computed $f_{P_n}(\x)$ for several sequences $\{P_n\}_{n \ge 1}$ of posets by developing MacMahon's
partition analysis. Corteel, Savage et. al. \cite{ Davis-05a, Davis-06}
presented five guidelines for deriving recurrence relations of $f_{P_n}(\x)$.
D'Souza \cite{Souza}
provided a Maple package {\tt GFPartitions} which generates
recurrence relations of $f_{P_n}(\x)$ once the decomposition of the posets $P_n$ is given manually.
Ekhad and Zeilberger \cite{Ekhad-07} considered the umbral operator
on ``grafting'' of posets. The corresponding Maple package {\tt
RotaStanley} can generate the recurrence relations automatically.

Our main objective is to find an efficient method of computing $f_P(\x)$.
For this purpose, we introduce two kinds of transformations on posets in Section \ref{sec-transform}, which are the \emph{deletion}
and the \emph{partially linear extension}. We find that there exist simple relations between the generating function of $P$-partitions of a poset and those of its transformations. Thus the generating function $f_P(\x)$ can be computed by recursively applying these transformations.

We then consider the posets $P_n$ composed of $n$ copies of a given poset $P$ in Section \ref{sec-gfun}. More precisely, we introduce the \emph{partially ordinal sum} $\oplus_R$ of posets, and denote by $P_n$ the sum $P \oplus_R P \oplus_R \dots \oplus_R P$, where $P$ occurs $n$ times. By applying the above two transformations, we find that the sequence $\{f_{P_n}(\x)\}_{n \ge 1}$  satisfies a system of recurrence relations with respect to $n$.

Finally, we provide some examples in Section \ref{sec-examples}, including a kind of $3$-rowed posets which can not be dealt with by the packages {\tt GFPartitions} and {\tt RotaStanley}.

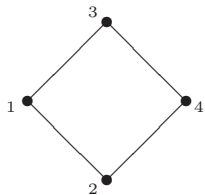
\begin{figure}[ht]
\centering
\begin{picture}(90,70)(0,0)

  \multiput(15,35)(30,-30){2}{\line(1,1){30}}
  \multiput(15,35)(30,30){2}{\line(1,-1){30}}

  \multiput(15,35)(60,0){2}{\circle*{4}}
  \multiput(45,5)(0,60){2}{\circle*{4}}

  \put(7,31){\tiny $1$}\put(38,0){\tiny $2$}
  \put(38,67){\tiny $3$}\put(78,31){\tiny $4$}

  \end{picture}
\caption{The Hasse diagram of the Diamond poset $D$. }\label{Fig-diamond}
\end{figure}

Before our further discussion, let us recall the Hasse diagram of a poset $(P,\leq)$. For $x,y\in P$, we say $y$
covers $x$, denoted by $x \lessdot y$, if $x < y$ and if no
element $z\in P$ satisfies $x<z<y$. The Hasse diagram of $P$ is a graphical representation of $P$, in which every element of
$P$ is represented by a vertex and two vertices are joined by a line
with vertex $x$ drawn below vertex $y$ if $x \lessdot y$.
To coincide with
the descriptions used by Andrews, Paule and Riese \cite{Andre01-a},
we rotate the Hasse diagram by $90$ degree clockwise so that smaller elements lie to the left. For
example, the diamond poset $D=\{1,2,3,4\}$ with cover relations $\{1 \lessdot  2 \lessdot 4, 1 \lessdot 3 \lessdot 4\}$ can be represented by Figure \ref{Fig-diamond}.

\section{Two Transformations on Posets}
\label{sec-transform}

In this section, we introduce the deletions and the partially linear extensions. On one hand, they reduce a poset to simpler ones. On the other hand, the generating function of $P$-partitions of a poset can be expressed in terms of those of its transformations.
This enables us to compute $f_P(\x)$ by recursively transforming the poset $P$.

\subsection{The Deletions}
The first transformation we consider is removing an element from a poset. We only concern on the {\it removable} elements, which cover at
most one element and are covered by at most one element.

Let $P$ be a poset and $b\in P$ be a removable element of $P$. The deletion of $P$ with respect to $b$ is the transformation from $P$ to the induced sub-poset $P \setminus \{b\}$ of $P$ by deleting $b$ from $P$. The following theorem shows the relation between $f_P(\x)$ and $f_{P \setminus \{b\}}(\x)$.

\begin{theo}\label{theo-absorb}
Let $P$ be a poset and $b \in P$ be a removable element of $P$. Then
\begin{equation}\label{eq-absorb}
f_P(\x) = \frac{g(\x) - h(\x)}{1-x_b},
\end{equation}
where
\[
g(\x) = \begin{cases}
f_{P \setminus \{b\}}(\x)|_{x_c=x_b x_c}, & \mbox{if $\exists\, c \in P$ such that $b \lessdot c$}, \\[5pt]
f_{P \setminus \{b\}}(\x), & \mbox{otherwise,}
\end{cases}
\]
and
\[
h(\x) = \begin{cases}
x_b f_{P \setminus \{b\}}(\x)|_{x_a=x_a x_b}, & \mbox{if $\exists\, a \in P$ such that $a \lessdot b$}, \\[5pt]
0, & \mbox{otherwise.}
\end{cases}
\]
\end{theo}
\noindent {\it Proof.} We only give the proof for the case when there exist $a$ and $c$ such that $a \lessdot b
\lessdot c$. The other three cases can be proved in a similar way.

By the definition of $P$-partitions, we have
\begin{eqnarray*}
f_P(\x)
&=& \sum_{\sigma \in Par(P \setminus \{b\})} \left( \sum_{m=\sigma(c)}^{\sigma(a)} x_b^m \prod_{u \in P \setminus \{b\}} x_u^{\sigma(u)}
\right) \\
&=& \sum_{\sigma \in Par(P \setminus \{b\})} \left( \frac{x_b^{\sigma(c)} - x_b^{\sigma(a)+1}}{1-x_b} \prod_{u \in P \setminus \{b\}} x_u^{\sigma(u)}
\right) \\
&=& \frac{1}{1-x_b} \sum_{\sigma \in Par(P \setminus \{b\})} (x_b
x_c)^{\sigma(c)} \prod_{u \in P \setminus \{b,c\}}
x_u^{\sigma(u)}\\
&& \quad - \frac{x_b}{1-x_b} \sum_{\sigma \in Par(P \setminus \{b\})}
(x_a x_b)^{\sigma(a)} \prod_{u \in P \setminus \{a,b\}}
x_u^{\sigma(u)},
\end{eqnarray*}
as desired. \qed

Now we give an example to illustrate the usage of Theorem~\ref{theo-absorb}.
\begin{exam}\label{exam-delete-gfun}
Let us consider the Diamond poset $D$ as shown in Figure~\ref{Fig-diamond}. We see that $2$ is a removable element of $D$ and $D \setminus \{2\}$ is a chain $C$ of length $3$ with
\[
f_C(x_1,x_2,x_3) = \frac{1}{(1-x_1)(1-x_1x_2)(1-x_1x_2x_3)}.
\]
Invoking equation \eqref{eq-absorb}, we derive that
\begin{align*}\label{eq-exam-absorb}
f_{D}(\x) &= \frac{1}{1-x_2}
  \left(f_C(x_1,x_3,x_2x_4)-x_2 f_C(x_1x_2,x_3,x_4)\right)\\[5pt]
 &= \frac{1-x_1^2x_2x_3}{(1-x_1)(1-x_1x_2)(1-x_1x_3)(1-x_1x_2x_3)(1-x_1x_2x_3x_4)}.
\end{align*}

\end{exam}

\subsection{The Partially Linear Extensions}
The second transformation we consider is partially ordering the elements of an anti-chain of a poset.

Let $A$ be an anti-chain of a poset $P$ and let $M$ be a non-empty subset of $A$.
The partially linear extension (PLE in short) of $P$ with respect to the pair $(M,A)$ is the transformation from the poset $P$ to the poset $P(M,A)$ by gluing the elements of $M$ together and
setting the glued element cover the elements of $A\setminus M$.
More precisely, $P(M,A)$ is the poset
defined on $P \setminus M \cup \{M\}$ and partially
ordered by $x \le y$ if and only if
\begin{itemize}
\item[(a)] $x,y \in P\setminus M$ and $x \le_P y$, or
\item[(b)] $x,y \in P\setminus M$ and there exist $x' \in A$ and $y' \in M$
     such that $x \le_P x', y' \le_P y$, or
\item[(c)] $x=M$ and there exists $y' \in M$ such that $y' \le_P y$,
or
\item[(d)] $y=M$ and there exists $x' \in A$ such that $x \le_P x'$,
or
\item[(e)] $x=y=M$.
\end{itemize}

\begin{exam} \label{exam-PLE}
Let $P=\{1,2,3,4,5\}$ be the poset as shown in Figure \ref{Fig-PLE0}.
\begin{figure}[ht]
\centering
\begin{picture}(70,70)(0,0)

  \multiput(5,35)(30,-30){2}{\line(1,1){30}}
  \multiput(5,35)(30,30){2}{\line(1,-1){30}}
  \put(5,35){\line(1,0){60}}

  \multiput(5,35)(30,0){3}{\circle*{4}}
  \multiput(35,5)(0,60){2}{\circle*{4}}

  \put(0,30){\tiny $1$}\put(37,0){\tiny $2$}
  \put(37,28){\tiny $3$}
  \put(37,67){\tiny $4$}\put(67,30){\tiny $5$}

  \end{picture}
\caption{The Hasse diagram of the poset $P$. }\label{Fig-PLE0}
\end{figure}
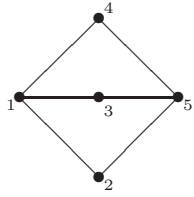

Then the posets
$P(\{2\},\{2,3,4\})$ and $P(\{2,3 \},\{2,3,4\})$ are given as in Figure~\ref{Fig-ple}.
\end{exam}
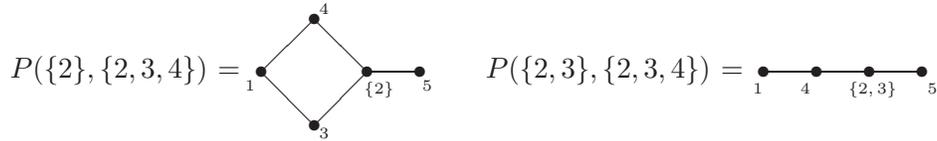
\begin{figure}[ht]
\centering
\begin{picture}(350,50)(0,0)

\put(0,23){\small$P(\{2\},\{2,3,4\})=$}
\multiput(95,25)(40,0){2}{\circle*{4}}
\multiput(115,5)(0,40){2}{\circle*{4}}
\multiput(115,5)(-20,20){2}{\line(1,1){20}}
\multiput(115,5)(20,20){2}{\line(-1,1){20}}
\put(155,25){\circle*{4}}
\put(135,25){\line(1,0){20}}
\put(89,18){\tiny$1$}
\put(117,0){\tiny$3$}
\put(117,47){\tiny$4$}
\put(134,17){\tiny$\{2\}$}
\put(156,18){\tiny$5$}

\put(180,23){\small$P(\{2,3\},\{2,3,4\})=$}
\put(285,25){\line(1,0){60}}
\multiput(285,25)(20,0){4}{\circle*{4}}
\put(281,17){\tiny$1$}
\put(299,17){\tiny$4$}
\put(317,17){\tiny$\{2,3\}$}
\put(347,17){\tiny$5$}
\end{picture}
\caption{Two partially linear extensions of $P$. \label{Fig-ple}}
\end{figure}

As a generalization of Theorem~2.3 in \cite{GHX-07}, the generating
function of $P$-partitions of a poset can be expressed by those
of its PLE's.
\begin{theo}\label{theo-PLE}
Let $P$ be a poset and $A$ be an anti-chain of $P$. Then
\begin{equation}\label{eq-PLE}
f_P(\x) = \sum_{\emptyset \not= M \subseteq A} (-1)^{|M|-1}
f_{P(M,A)}(\x) |_{x_M = \prod_{a \in M}x_a}.
\end{equation}
\end{theo}
\noindent {\it Proof.} For each $a \in A$, we define
\[
S_a = \{\sigma \in \pa(P) \colon \sigma(a) \le \sigma(x), \forall\,
x \in A \}.
\]
Since $\pa(P) = \cup_{a \in A} S_a$, by the
inclusion-exclusion principle we derive that
\[
f_P(\x) = \sum_{\emptyset \not= M \subseteq A} (-1)^{|M|-1}
\sum_{\sigma \in S_M} \prod_{u \in P} x_u^{\sigma{(u)}},
\]
where $S_M = \bigcap_{a \in M} S_a$.

Given a $P$-partition $\sigma \in S_M$, we denote $m=\min\{\sigma(a)
\colon a \in A\}$. Then by the definition of $S_M$ we have
$\sigma(x)=m$ for any $x \in M$. Let
\[
\sigma'(u) = \begin{cases} \sigma(u), & u \in P \setminus M, \\
m, & u=M. \end{cases}
\]
One sees that $\sigma'$ is a $P$-partition of $P(M,A)$.
Conversely, let $\sigma'$ be a $P$-partition of $P(M,A)$. By
defining
\[
\sigma(u) = \begin{cases} \sigma'(u), & u \in P \setminus M, \\
\sigma'(M), & u \in M, \end{cases}
\]
we obtain a $P$-partition of $P$ in $S_M$. We thus set up a one-to-one
corresponding between the $P$-partitions of $P$ in $S_M$ and the
$P$-partitions of $P(M,A)$. Therefore,
\begin{eqnarray*}
f_P(\x) &=& \sum_{\emptyset \not= M \subseteq A} (-1)^{|M|-1}
\sum_{\sigma' \in \pa(P(M,A))} \prod_{u \in M} x_u^{\sigma'(M)}
\cdot \prod_{u \in P \setminus M} x_u^{\sigma'{(u)}} \\
&=& \sum_{\emptyset \not= M \subseteq A} (-1)^{|M|-1} f_{P(M,A)}(\x) |_{x_M = \prod_{a \in M}x_a}.
\end{eqnarray*}
This completes the proof. \qed

The Hasse diagrams provide a simple graphical representation for Equation~\eqref{eq-PLE}.
\begin{exam}
Let $P=\{1,2,3,4\}$ be the poset as shown in Figure \ref{Fig-exam-PLE}.
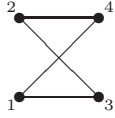
\begin{figure}[htb]
\centering
\begin{picture}(60,40)(0,0)

  \multiput(5,5)(0,30){2}{\line(1,0){30}}
  \multiput(5,5)(30,0){2}{\circle*{4}}
  \multiput(5,35)(30,0){2}{\circle*{4}}
  \put(5,5){\line(1,1){30}}
  \put(35,5){\line(-1,1){30}}
  \put(0,0){\tiny$1$}
  \put(0,37){\tiny$2$}
  \put(37,0){\tiny$3$}
  \put(37,37){\tiny$4$}
 \end{picture}
\caption{The Hasse diagram of the poset $P$.}\label{Fig-exam-PLE}
\end{figure}
Taking the anti-chain $\{1,2\}$ into account, we find that the PLE's are shown as in Figure \ref{Fig-exam-PLE-gfun}, from which we read out
\begin{multline*}
f_P(\x) =  f_{P(\{1\},\{1,2\})}(x_2,x_1,x_3,x_4) +  f_{P(\{2\},\{1,2\})}(x_1,x_2,x_3,x_4) \\
           -f_{P(\{1,2\},\{1,2\})}(x_1x_2, x_3, x_4).
\end{multline*}

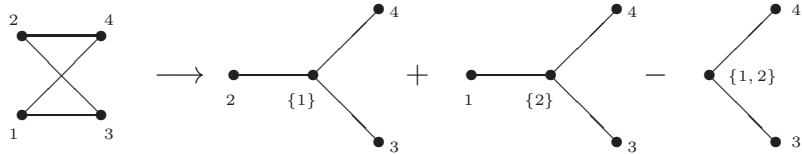
\begin{figure}[ht]
\centering
\begin{picture}(310,60)(0,0)

  \multiput(5,15)(0,30){2}{\line(1,0){30}}
  \multiput(5,15)(30,0){2}{\circle*{4}}
  \multiput(5,45)(30,0){2}{\circle*{4}}
  \put(5,15){\line(1,1){30}}
  \put(35,15){\line(-1,1){30}}
  \put(0,6){\tiny$1$}
  \put(0,49){\tiny$2$}
  \put(36,6){\tiny$3$}
  \put(36,49){\tiny$4$}

  \put(55,27){$\longrightarrow$}

  \multiput(85,30)(30,0){2}{\circle*{4}}
  \put(85,30){\line(1,0){30}}
  \put(115,30){\line(1,1){25}}
  \put(115,30){\line(1,-1){25}}
  \multiput(140,5)(0,50){2}{\circle*{4}}
  \put(82,19){\tiny$2$}
  \put(105,19){\tiny$\{1\}$}
  \put(144,1){\tiny$3$}
  \put(144,52){\tiny$4$}

  \put(150,27){$+$}

  \multiput(175,30)(30,0){2}{\circle*{4}}
  \put(175,30){\line(1,0){30}}
  \put(205,30){\line(1,1){25}}
  \put(205,30){\line(1,-1){25}}
  \multiput(230,5)(0,50){2}{\circle*{4}}
  \put(172,19){\tiny$1$}
  \put(195,19){\tiny$\{2\}$}
  \put(234,1){\tiny$3$}
  \put(234,52){\tiny$4$}

  \put(240,27){$-$}

  \put(265,30){\circle*{4}}
  \put(265,30){\line(1,1){25}}
  \put(265,30){\line(1,-1){25}}
  \multiput(290,5)(0,50){2}{\circle*{4}}
  \put(272,28){\tiny$\{1,2\}$}
  \put(296,3){\tiny$3$}
  \put(296,53){\tiny$4$}
 \end{picture}
\caption{A graphical representation of Equation~\eqref{eq-PLE}.}\label{Fig-exam-PLE-gfun}
\end{figure}
\end{exam}

\subsection{Computing $f_P(\x)$ via two Transformations}
We shall show that the deletion and the partially linear extension are powerful enough for computing $f_P(\x)$ for any poset $P$.

\begin{theo}\label{theo-transformations}
Any poset can be reduced to the empty poset by applying the deletion and the partially linear extension finite times.
\end{theo}
\pf Let $ac(P)$ denote the number of distinct anti-chains of a poset $P$. Since each element of $P$ forms an anti-chain, the deletion reduces $ac(P)$ by at least one.

Suppose that $A$ is an anti-chain of $P$ with cardinality at least two and $M$ is a non-empty subset of $A$. If $M \not= A$, there exist $x \in M$ and $y \in A \setminus M$. Then $\{x,y\}$ is an anti-chain of $P$ but is not an anti-chain of $P(M,A)$. If $M = A$, the cardinality of $P(M,A)$ is strictly less than that of $P$. Thus in either case, we have $ac(P(M,A)) \le ac(P)-1$.

Now iteratively apply the deletion whenever there is a removable element and apply the partially linear extension whenever there is an anti-chain with cardinality at least two. Since $ac(P)$ is a finite number, the procedure eventually stops. The final poset contains no removable element and no anti-chain with cardinality at least two. The only poset satisfying this property is the empty poset. \qed

\section{Partially Ordinal Sums}\label{sec-gfun}

In this section, we consider a kind of posets composed of small blocks.

Let $P,Q$ be two posets and $R$ be a subset of the Cartesian product
$P\times Q$. The \emph{partially ordinal sum} (or \emph{$R$-plus}, for short) of $P$ and $Q$ with respect to $R$
is the poset $P \oplus_R Q$ defined on the disjoint union of $P$ and
$Q$ and partially ordered by $x \leq y$ in $P \oplus_R Q$ if and
only if
\begin{itemize}
\item[(a)] $x,y \in P$ and $x \leq_P y$, or
\item[(b)] $x,y\in Q$ and $x \leq_Q y$, or
\item[(c)] $x \in P, y \in Q$ and there exists $(x',y') \in R$ such that $x \le_P x'$ and
$y' \le_Q y$.
\end{itemize}
As special cases, $R$-plus reduces to the direct sum if $R=\emptyset$ and to the ordinal sum if $R=P \times Q$, respectively.

\begin{exam}\label{exam-R-plus}
Figure~\ref{R-plus} gives the $R$-plus of $P$ and $Q$, where $P=\{1,2\}$ is an anti-chain, $Q=\{3,4\}$ is a chain, and $R=\{(1,4), (2,3), (2,4)\}$.
\end{exam}
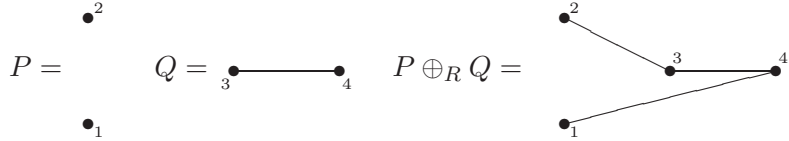
\begin{figure}[ht]
\centering
\begin{picture}(300,50)(0,0)

  \put(0,23){\small$P=$}
  \multiput(30,5)(0,40){2}{\circle*{4}}
  \put(32,0){\tiny$1$}
  \put(32,47){\tiny$2$}

  \put(55,23){\small$Q=$}
  \put(85,25){\line(1,0){40}}
  \multiput(85,25)(40,0){2}{\circle*{4}}
  \put(80,18){\tiny$3$}
  \put(126,18){\tiny$4$}

  \put(145,23){\small$P \oplus_R Q =$}
  \multiput(210,5)(0,40){2}{\circle*{4}}
  \multiput(250,25)(40,0){2}{\circle*{4}}
  \put(210,5){\line(4,1){80}}
  \put(210,45){\line(2,-1){40}}
  \put(250,25){\line(1,0){40}}
  \put(212,0){\tiny$1$}
  \put(212,47){\tiny$2$}
  \put(251,28){\tiny$3$}
  \put(291,28){\tiny$4$}
\end{picture}
\caption{The partially ordinal sum of $P$ and $Q$ w.r.t.
$R$. \label{R-plus}}
\end{figure}

It is easy to check that the partially ordinal sum is associative up to  isomorphic. Therefore we can naturally extend the definition of partially ordinal sum of two posets to several posets.
In particular,  we denote by $P_R^n$ the partially ordinal sum $P \oplus_R P \oplus_R \cdots \oplus_R P$ of
$n$ copies of $P$ w.r.t. $R$.

The sequence $\{f_{P_R^n}(\x)\}_{n \ge 1}$ satisfies a kind of recurrence relation which is given as follows. We say a sequence $\{f_n(\x)\}_{n\ge1}$ of functions is \emph{substituted recursive} if there are finitely many sequences
\[
\{f_n^{(0)}(\x)\}_{n \ge 1},\ \{f_n^{(1)}(\x)\}_{n \ge 1},\ \ldots,\ \{f_n^{(I)}(\x)\}_{n \ge 1}
\]
such that
$f_n^{(0)}(\x) = f_n(\x)$ and for $i=0,1,\ldots,I$,
\begin{equation} \label{eq-subrec}
 f_n^{(i)}(\x) = \sum _{j=0}^I \sum _{k=0}^{K} r_{i,j,k}(\x) f_{n-1}^{(j)}(\y^{(j,k)}),
\end{equation}
where $r_{i,j,k}(\x)$ are rational functions and each component of the variable vector $\y^{(j,k)}$ is a  monomial in $x_1,\ldots,x_n$.

For example, suppose that
\[
f_n(x_1,\ldots,x_n) = \frac{1}{1-x_1} g_{n-1}(x_2, \ldots,x_n),
\]
and
\[
g_n(x_1,\ldots,x_n) = \frac{g_{n-1}(x_1x_2, x_3, \ldots,x_n)}{1-x_1x_2} - \frac{g_{n-1}(x_1, x_2x_3, x_4, \ldots, x_n)}{1-x_1x_3}.
\]
Then both $\{f_n(\x)\}_{n \ge 1}$ and $\{g_n(\x)\}_{n \ge 1}$ are substituted recursive.

To compute $f_{P_R^n}(\x)$, we consider the more general posets
\begin{equation}\label{eq-Xn}
 X_n= A \oplus_{R_1} P_{R}^n \oplus_{R_2} B,
\end{equation}
where $A,B$ are posets and $R_1 \subseteq A \times P, R_2 \subseteq P \times B$.

\begin{theo}
\label{theo-subrec}
Let $X_n$ be given as in \eqref{eq-Xn}.
Then the sequence $\{f_{X_n}(\x)\}_{n \ge 1}$ of generating functions of $P$-partitions of
$X_n$ is substituted recursive.
\end{theo}

\pf
Let ${\cal S}$ denote the set of all pairs $(C,R')$ such that $C$ is a chain, $R' \subseteq C\times P$, and none of the elements of $C$ in $C \oplus_{R'} P$ is removable. Since each element of $C$ is not removable, it must be covered by a certain element of $P$. Moreover, two distinct elements of $C$ can not be covered by the same element in $P$. This implies that the cardinality of $C$ is less than that of $P$. Therefore, ${\cal S}$ is a finite set.

Now we apply the deletion and the partially linear extension to the elements of $A \oplus_{R_1} P$ in $X_n$ whenever possible. As shown in the proof of Theorem~\ref{theo-transformations}, we eventually arrive at posets of the form $C \oplus_{R'} P_R^{n-1} \oplus_{R_2} B$ with $(C,R') \in {\cal S}$. Moreover, we have
\[
f_{X_n}(\x) = \sum_{(C,R') \in \cal S} r(C,R',\x) f_{C \oplus_{R'} P_R^{n-1} \oplus_{R_2} B}(\y^{(C,R')}),
\]
where $r(C,R',\x)$ are rational functions of $\x$ depending on $C$ and $R'$ and each component of the variable vector $\y^{(C,R')}$ is a monomial in $\x$. By a similar discussion, for each $(C,R') \in {\cal S}$ we have
\[
f_{C \oplus_{R'} P_R^{n} \oplus_{R_2} B}(\x) = \sum_{(C',R'') \in \cal S} r'(C,R', C',R'',\x) f_{C' \oplus_{R''} P_R^{n-1} \oplus_{R_2} B}(\y^{(C,R',C',R'')}),
\]
where $r'(C,R',C',R'',\x)$ are rational functions of $\x$ that depend on $C,R',C'$ and $R''$. This completes the proof. \qed

Note that the proof of Theorem \ref{theo-subrec} provides an algorithm for generating substituted
recurrence relations of $f_{X_n}(\x)$. The corresponding {\tt Maple} package is available at
{\tt http://www.combinatorics.net.cn/homepage/hou/}.

\section{Some Examples}\label{sec-examples}

In this section, we present some examples to illustrate our approach to the computation of $f_{X_n}(\x)$.
We begin with an introductory example, i.e., the $3$-rowed plane partition introduced by Souza \cite{Souza}. Then we provide some more examples, including the zigzag posets, the 2-rowed posets with double diagonals and the multi-cube posets.

In these examples, we consider the $q$-generating function $f_P(q)$ obtained from $f_P(\x)$ by setting all variables $x_i$ equal the indeterminant $q$. For brevity, we omit some variables equalling $q$ and write $f(x_1,x_2,\ldots,x_k)$ instead of $f(x_1,x_2,\ldots,x_k,q,q,\ldots,q)$. We also adopt the standard notation
\begin{eqnarray*}
(a;q)_n &:=& (1-a)(1-aq)\cdots(1-aq^{n-1}),\\
(a;q)_\infty &:=& \prod_{n=0}^\infty (1-aq^n).
\end{eqnarray*}

\subsection{An Introductory Example}

\begin{figure}[ht]
\centering
\begin{picture}(170,85)(0,0)

  \multiput(10,40)(30,0){3}{\circle*{4}}
  \multiput(40,10)(30,0){3}{\circle*{4}}
  \multiput(40,70)(30,0){3}{\circle*{4}}

  \multiput(10,40)(30,0){3}{\line(1,1){30}}
  \multiput(10,40)(30,0){3}{\line(1,-1){30}}

  \multiput(40,10)(0,60){2}{\line(1,0){70}}

  \multiput(105,40)(6,0){3}{\circle*{2}}

  \put(130,40){\circle*{4}}
  \put(130,40){\line(1,1){30}}
  \put(130,40){\line(1,-1){30}}
  \multiput(160,10)(0,60){2}{\circle*{4}}

  \put(160,10){\line(-1,0){10}}
  \put(160,70){\line(-1,0){10}}

  \put(2,37){\tiny$1$}
  \put(39,0){\tiny$2$}
  \put(39,74){\tiny$3$}
  \put(32,37){\tiny$4$}
  \put(69,0){\tiny$5$}
  \put(69,74){\tiny$6$}
  \put(62,37){\tiny$7$}
  \put(97,0){\tiny$8$}
  \put(99,74){\tiny$9$}
  \put(135,38){\tiny$3n-2$}
  \put(148,2){\tiny$3n-1$}
  \put(158,74){\tiny$3n$}

\end{picture}
\caption{The graphical representation of the $3$-rowed poset $P_n$.}\label{Fig-exam-3-rowed}
\end{figure}
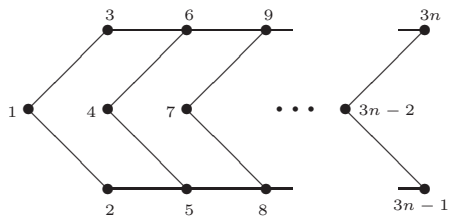

D'Souza\cite{Souza} introduced the $3$-rowed plane partition whose corresponding poset $P_n$ is given by Figure \ref{Fig-exam-3-rowed}.
He failed to find out recurrence relations of the generating function $f_{P_n}(\x)$.
Our approach gives substituted recurrence relations of $f_{P_n}(\x)$.

It is easy to see that $P_n = P_R^n$, where
$P$ is the poset on $\{1,2,3\}$ with $1 \lessdot 2$ and $1 \lessdot 3$,  and $R=\{(2,2),(3,3)\}$.

\begin{figure}[ht]
\centering
\begin{picture}(270,50)(0,0)

  \multiput(5,25)(20,0){3}{\circle*{2}}
  \put(25,5){\circle{3}}
  \multiput(45,5)(20,0){2}{\circle*{2}}
  \multiput(25,45)(20,0){3}{\circle*{2}}

  \multiput(5,25)(20,0){3}{\line(1,1){20}}
  \multiput(5,25)(20,0){3}{\line(1,-1){20}}

  \multiput(25,5)(0,40){2}{\line(1,0){47}}

  \multiput(60,25)(6,0){3}{\circle*{1}}

  %\put(80,25){\tiny Absorbtion}
  \put(80,22){$\rightarrow$}

  \multiput(105,25)(20,0){3}{\circle*{2}}
  \multiput(145,5)(20,0){2}{\circle*{2}}
  \put(125,45){\circle{3}}
  \multiput(145,45)(20,0){2}{\circle*{2}}

  \multiput(105,25)(20,0){3}{\line(1,1){20}}
  \multiput(125,25)(20,0){2}{\line(1,-1){20}}

  \put(105,25){\line(2,-1){40}}

  \put(145,5){\line(1,0){27}}
  \put(125,45){\line(1,0){47}}

  \multiput(160,25)(6,0){3}{\circle*{1}}

  %\put(180,25){\tiny Absorbtion}
  \put(180,22){$\rightarrow$}

  \multiput(205,25)(20,0){3}{\circle*{2}}
  \multiput(245,5)(20,0){2}{\circle*{2}}
  \multiput(245,45)(20,0){2}{\circle*{2}}

  \multiput(225,25)(20,0){2}{\line(1,1){20}}
  \multiput(225,25)(20,0){2}{\line(1,-1){20}}

  \put(205,25){\line(2,-1){40}}
  \put(205,25){\line(2,1){40}}

  \put(245,5){\line(1,0){27}}
  \put(245,45){\line(1,0){27}}

  \multiput(260,25)(6,0){3}{\circle*{1}}

\end{picture}
\caption{The transformation from $P_n$ to $Q_{n-1}$.}\label{Fig-exam-3-rowed-absorb}
\end{figure}
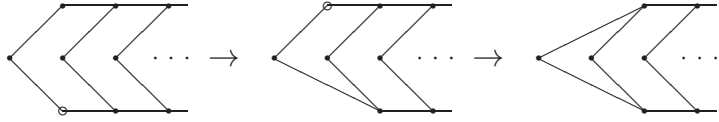

By deleting the removable elements $2$ and $3$ of
$P$, we reduce the poset $P_n$ to $
Q_{n-1} = {\bf 1} \oplus_{R'} P_R^{n-1}$,
where ${\bf 1}$ is the poset with only one element and $R'=\{(1,2),(1,3)\} \subset {\bf 1} \times P$. See Figure~\ref{Fig-exam-3-rowed-absorb} for a demonstration.
According to Theorem \ref{theo-absorb}, we derive that
\begin{multline}
f_{P_n}(x_1,x_2,x_3) = \frac{1}{(1-x_2)(1-x_3)}\\[7pt] \times \big(
f_{Q_{n-1}}(x_1,q,qx_2,qx_3)  -x_2 f_{Q_{n-1}}(x_1x_2,q,q,qx_3)\\
 -x_3 f_{Q_{n-1}}(x_1 x_3,q,qx_2,q)  + x_2x_3f_{Q_{n-1}}(x_1x_2 x_3,q,q,q) \big).
\label{eq-3-rowed-subrec1}
\end{multline}

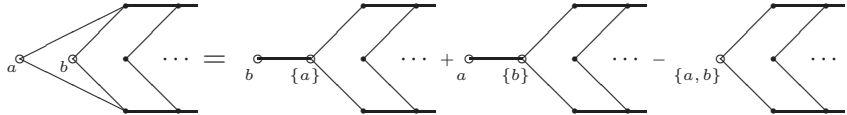
\begin{figure}[ht]
\centering
\begin{picture}(320,50)(0,0)

  \multiput(5,25)(20,0){2}{\circle{3}}
  \put(45,25){\circle*{2}}
  \multiput(45,5)(20,0){2}{\circle*{2}}
  \multiput(45,45)(20,0){2}{\circle*{2}}

  \multiput(25,25)(20,0){2}{\line(1,1){20}}
  \multiput(25,25)(20,0){2}{\line(1,-1){20}}
  \multiput(45,5)(0,40){2}{\line(1,0){27}}

  \put(5,25){\line(2,-1){40}}
  \put(5,25){\line(2,1){40}}

  \multiput(60,25)(4,0){3}{\circle*{1}}

  \put(0,20){\tiny$a$}
  \put(20,20){\tiny$b$}

  \put(74,22){$=$}

  \multiput(95,25)(20,0){2}{\circle{3}}
  \put(135,25){\circle*{2}}
  \multiput(135,5)(20,0){2}{\circle*{2}}
  \multiput(135,45)(20,0){2}{\circle*{2}}

  \multiput(115,25)(20,0){2}{\line(1,1){20}}
  \multiput(115,25)(20,0){2}{\line(1,-1){20}}
  \multiput(135,5)(0,40){2}{\line(1,0){27}}

  \put(95,25){\line(1,0){20}}

  \multiput(150,25)(4,0){3}{\circle*{1}}

  \put(90,17){\tiny$b$}
  \put(107,17){\tiny$\{a\}$}

  \put(160,23){ \tiny$+$}

  \multiput(175,25)(20,0){2}{\circle{3}}
  \put(215,25){\circle*{2}}
  \multiput(215,5)(20,0){2}{\circle*{2}}
  \multiput(215,45)(20,0){2}{\circle*{2}}

  \multiput(195,25)(20,0){2}{\line(1,1){20}}
  \multiput(195,25)(20,0){2}{\line(1,-1){20}}
  \multiput(215,5)(0,40){2}{\line(1,0){27}}

  \put(175,25){\line(1,0){20}}

  \multiput(230,25)(4,0){3}{\circle*{1}}

  \put(170,17){\tiny$a$}
  \put(187,17){\tiny$\{b\}$}

  \put(240,23){ \tiny$-$}

  \put(270,25){\circle{3}} \put(290,25){\circle*{2}}
  \multiput(290,45)(20,0){2}{\circle*{2}}
  \multiput(290,5)(20,0){2}{\circle*{2}}

  \multiput(270,25)(20,0){2}{\line(1,1){20}}
  \multiput(270,25)(20,0){2}{\line(1,-1){20}}
  \multiput(290,5)(0,40){2}{\line(1,0){27}}

  \multiput(305,25)(4,0){3}{\circle*{1}}

  \put(252,17){\tiny$\{a,b\}$}

\end{picture}
\caption{The PLE transformation of $Q_n$.}\label{Fig-exam-3-rowed-PLE}
\end{figure}

Now let us consider $Q_n={\bf 1} \oplus_{R'} P_R^n$. It is readily to see that the unique element of ${\bf 1}$ and the minimal element of $P$ are not comparable. Thus we can apply PLE to the anti-chain consisting of these two elements, as shown in Figure~\ref{Fig-exam-3-rowed-PLE}.

After further deletions, all the posets generated by the PLE transformation reduce to $Q_{n-1}$. We thus obtain the recurrence relation
\begin{multline}
f_{Q_n}(x_1,x_2,x_3,x_4) = \frac{1-x_1x_2}{(1-x_1)(1-x_2)(1-x_3)(1-x_4)} \\[7pt]
\times  \big( f_{Q_{n-1}}(x_1 x_2,q, qx_3,qx_4) -x_3f_{Q_{n-1}}(x_1x_2 x_3,q,q,qx_4) \\[7pt]
 -x_4f_{Q_{n-1}}(x_1 x_2 x_4, q,qx_3,q)+x_3x_4f_{Q_{n-1}}(x_1x_2 x_3x_4,q,q,q) \big).
\label{eq-3-rowed-subrec2}
\end{multline}

By the recurrence relations \eqref{eq-3-rowed-subrec1} and \eqref{eq-3-rowed-subrec2} and the initial condition
\[
 f_{Q_1}(x_1,x_2,x_3,x_4) = \frac{1-x_1^2x_2^2x_3x_4}{(1-x_1)(1-x_2)(1-x_1x_2x_3)(1-x_1x_2x_4)(1-x_1x_2x_3x_4)},
\]
we can compute $f_{P_n}(q)$ recursively.

\subsection{More Examples}
In this subsection, we give three more examples: the zigzag posets,
the $2$-rowed posets with double diagonals and the multi-cube posets.

\begin{exam}
Let $P=\{1,2\}$ be a chain with $2 \lessdot 1$ and let $R=\{(2,1)\}$.
The zig-zag poset of length $n$ is given by $Z_n = P_R^n$. We have
\[
f_{Z_n}(x_1,x_2)  = \frac{f_{Z_{n-1}}(qx_2,q)}{(1-x_1)(1-x_2)}
  -\frac{x_1 f_{Z_{n-1}}(qx_1x_2,q)}{(1-x_1)(1-x_1x_2)}.
\]
The initial condition is given by
\[
f_{Z_1}(x_1,x_2) = \frac{1}{(1-x_2)(1-x_1x_2)}.
\]
\end{exam}
Note that the $P$-partitions of $Z_n$ is exactly the up-down sequences defined by Carlitz \cite{Carlitz1973}.

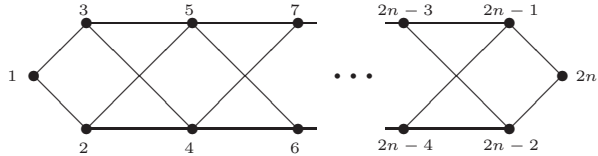
\begin{figure}[ht]
\centering
\begin{picture}(225,60)(0,0)

  \put(10,30){\circle*{4}}
  \put(10,30){\line(1,-1){20}}
  \put(10,30){\line(1,1){20}}
  \multiput(30,10)(40,0){3}{\circle*{4}}
  \multiput(30,50)(40,0){3}{\circle*{4}}
  \multiput(30,10)(0,40){2}{\line(1,0){87}}
  \multiput(30,10)(40,0){2}{\line(1,1){40}}
  \multiput(70,10)(40,0){2}{\line(-1,1){40}}

  \multiput(125,30)(6,0){3}{\circle*{2}}

  \multiput(150,10)(40,0){2}{\circle*{4}}
  \multiput(150,50)(40,0){2}{\circle*{4}}
  \multiput(143,10)(0,40){2}{\line(1,0){47}}
  \put(150,10){\line(1,1){40}}
  \put(190,10){\line(-1,1){40}}
  \put(190,10){\line(1,1){20}}
  \put(210,30){\circle*{4}}
  \put(210,30){\line(-1,1){20}}

  \put(0,28){\tiny$1$}
  \put(27,1){\tiny$2$}
  \put(27,53){\tiny$3$}
  \put(67,1){\tiny$4$}
  \put(67,53){\tiny$5$}
  \put(107,1){\tiny$6$}
  \put(107,53){\tiny$7$}
  \put(140,2){\tiny$2n-4$}
  \put(140,53){\tiny$2n-3$}
  \put(180,2){\tiny$2n-2$}
  \put(180,53){\tiny$2n-1$}
  \put(215,28){\tiny$2n$}

\end{picture}
\caption{The $2$-rowed poset with double diagonals.}\label{Fig-exam-2-rowed-doublediagonal}
\end{figure}

\begin{exam}\label{exam-2-rowed-doublediagonal}
Let $P_n$ be  the $2$-rowed poset with double diagonals depicted in Figure \ref{Fig-exam-2-rowed-doublediagonal}.
Then we have
\begin{equation}\label{eq-2-rowed}
f_{P_n}(q) = \frac{(-q^2;q^2)_{n-1}}{(q;q)_{2n}}.
\end{equation}
\end{exam}

Davis, Souza, Lee and Savage \cite{Davis-06} used the
``digraph method'' to derive formulae \eqref{eq-2-rowed}. Using the inclusion-exclusion
principle, Gao, Hou and Xin \cite{GHX-07} obtained the same generating function. Our approach leads to a recurrence relation as follows
\begin{eqnarray*}
&&f_{P_n}(x_1,x_2,x_3,\ldots,x_{2n})\\[5pt]
=&&\frac{1-x_1^2x_2x_3}{(1-x_1)(1-x_1x_2)(1-x_1x_3)}
    f_{P_{n-1}}(x_1x_2x_3,x_4,x_5,\ldots,x_{2n}).
\end{eqnarray*}
Note further that $P_2$ is the Diamond poset
$D$ given in Figure~\ref{Fig-diamond}. By iterating the recurrence relation and substituting all $x_i$ with $q$, we arrive at \eqref{eq-2-rowed}.

\begin{figure}[ht]
\centering
\begin{picture}(195,80)(0,0)

  \multiput(30,10)(60,0){3}{\circle*{4}}
  \multiput(50,30)(60,0){3}{\circle*{4}}
  \multiput(10,50)(60,0){3}{\circle*{4}}
  \multiput(30,70)(60,0){3}{\circle*{4}}
  \multiput(30,10)(20,20){2}{\line(1,0){140}}
  \multiput(10,50)(20,20){2}{\line(1,0){140}}
  \multiput(30,10)(60,0){3}{\line(1,1){20}}
  \multiput(10,50)(60,0){3}{\line(1,1){20}}
  \multiput(30,10)(60,0){3}{\line(-1,2){20}}
  \multiput(50,30)(60,0){3}{\line(-1,2){20}}

  \multiput(180,50)(6,0){3}{\circle*{2}}

  \put(3,42){\tiny$1$}
  \put(32,2){\tiny$2$}
  \put(23,73){\tiny$3$}
  \put(51,22){\tiny$4$}
  \put(63,42){\tiny$5$}
  \put(92,2){\tiny$6$}
  \put(83,73){\tiny$7$}
  \put(111,22){\tiny$8$}
  \put(123,42){\tiny$9$}
  \put(150,2){\tiny$10$}
  \put(143,73){\tiny$11$}
  \put(171,22){\tiny$12$}

\end{picture}
\caption{The graphical representation of multi-cube posets.}\label{Fig-exam-Cubics}
\end{figure}
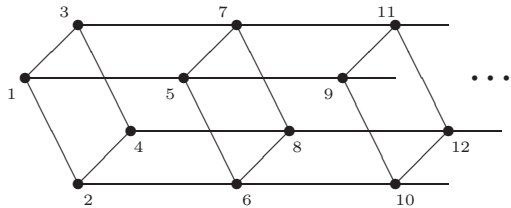

\begin{exam}
Let $D=\{1,2,3,4\}$ be the Diamond poset shown as in Figure~\ref{Fig-diamond} and
\[
R=\{(1,1), (2,2), (3,3), (4,4) \},
\]
The $n$-th multi-cube poset $C_n$ is defined by $C_n = D_R^n$, as shown in Figure \ref{Fig-exam-Cubics}.
Using the substituted recurrence relations, we compute $f_{C_n}(q)$ for $n \le 6$.
For example,
\[
f_{C_6}(q) = \frac{q^{192}+2q^{190}+\cdots+40660110q^{96}+\cdots+2q^2+1}{(q;q)_{24}}.
\]
\end{exam}

MacMahon \cite[Section $7$]{MacMahon-12} first considered the $P$-partitions of the poset $C_2$.
Under the framework of partition analysis, Andrews, Paule and Rises \cite{Andre01-c} computed the
$f_{C_2}(q)$ using the {\tt Omega} package. When $n \ge 3$, the computation of $f_{C_n}(q)$ seems not feasible by using {\tt Omega}.

\vskip 15pt

\noindent\textbf{Acknowledgments.}
This work was partially done during the second author's visiting at
 RISC-Linz. He would like to thank Peter Paule for his invitation.
This work was supported by the PCSIRT project of the Ministry of
Education and the National Science Foundation of China.


\begin{thebibliography}{99}


\bibitem{Andre01-c}
G.~E.~Andrews, P.~Paule and A.~Riese, MacMahon's partition analysis
III: the Omega package, Europ. J. Combin., {\bf 22} (2001) 887--904.

\bibitem{Andre01-a}
G.~E.~Andrews, P.~Paule and A.~Riese, MacMahon's partition analysis
VIII: plane partition diamonds, Adv. in Appl. Math., {\bf 27} (2001)
231--242.

\bibitem{Andre01-e}
G.~E.~Andrews, P.~Paule and A.~Riese, MacMahon's partition analysis
IX: $k$-gon partitions, Bull. Austral. Math. Soc., {\bf 64} (2001)
321--329.

\bibitem{Andre04-a}
G.~E.~Andrews, P.~Paule and A.~Riese, MacMahon's partition analysis
X: plane partitions with diagonals, South East Asian J. Math. Sci.,
{\bf 3} (2004) 3--14.

\bibitem{Andre04}
G.~E.~Andrews, P.~Paule and A.~Riese, MacMahon's partition analysis
XI: hexagonal plane partitions, 2004. SFB Report n. 2004--4, J.
Kepler University, Linz.

\bibitem{Andre07b}
G.~E.~Andrews and P.~Paule, MacMahon's partition analysis XII: plane
partitions, J. London Math. Soc., {\bf 76} (2007) 647--666.

\bibitem{Carlitz1973}
L.~Carlitz, Enumeration of up-down sequences. Discrete Math.,
{\bf 4} (1973) 273--386.

\bibitem{Davis-05a}
S.~Corteel, S.~Lee and C.~D.~Savage, Five guidelines for partition
analysis with applications to lecture hall-type theorems, Comb.
Number theory, {\bf 297} (2007) 131--155.

\bibitem{Davis-06}
J.~W.~Davis, E.~D'~Souza, S.~Lee and C.~D.~Savage, Enumeration of
integer solutions to linear inequalities defined by digraphs,
Integer points in polyhedra-geometry, number theory, representati
on theory, algebra, optimization, statistics, 79--91,
Contemp. Math., 452, Amer. Math. Soc., Providence, RI, 2008.

\bibitem{Ekhad-07}
S.~B.~Ekhad and D.~Zeilberger, Using Rota's umbral calculus to
enumerate to Stanley's $P$-partitions,
Adv. in Appl. Math. {\bf 41} (2008) 206-217.


\bibitem{GHX-07}
W.~Gao, Q.H.~Hou, and G.C.~Xin, On $P$-partitions related to
ordinal sums of posets, European J. Combin. {\bf 30 } (2009) 1370--1381.

\bibitem{MacMahon-12}
P.A.~MacMahon, Memoir on the theory of the partition of numbers --- Part VI,
 Phil. Trans., \textbf{211} (1912), 245--373.

\bibitem{Souza}
E.~D'~Souza,  Automating the enumeration of sequences defined by
digraphs, Thesis.

\bibitem{Stanley-97}
R.~P.~Stanley, Enumerative Combinatorics Vol. I, Cambridge
University Press, Cambridge, 1997.

\end{thebibliography}
\end{document}